\documentclass[a4paper]{article}

\usepackage{amsmath}
\usepackage{amssymb}
\usepackage{amsfonts}
\usepackage{enumerate}
\usepackage{vmargin}
\usepackage{verbatim}

\newcommand{\tens}{\otimes}

\newcommand{\M}{\ensuremath{\mathbb{M}}}

\newcommand{\spn}{\ensuremath{\mathop{\rm Span\,}\nolimits}}

\newcommand{\Id}{{\rm Id}}

\renewcommand{\leq}{\ensuremath{\leqslant}}
\renewcommand{\geq}{\ensuremath{\geqslant}}
\newcommand{\n}{\noindent}
\newcommand{\m}{\mathcal}
\newcommand{\qed}{\nobreak \hfill \vrule height6pt  width6pt depth0pt}

\newcommand{\cM}{\mathcal{M}}
\newcommand{\cN}{\mathcal{N}}
\newcommand{\ot}{\otimes}
\newcommand{\tr}{\tau}

\newtheorem{thm}{Theorem}[section]

\newtheorem{cor}[thm]{Corollary}
\newtheorem{lemma}[thm]{Lemma}

\newenvironment{rk}{{\noindent \bf Remark \addtocounter{thm}{1}\arabic{section}.\arabic{thm}\,}}
{
\smallskip
}
\newenvironment{pf}[1][]{\noindent {\it Proof #1} : }{\hbox{~}\qed
\smallskip
}

\title{Some remarks on noncommutative Khintchine inequalities}
\date{}
\author{Sjoerd Dirksen \and \'Eric Ricard}

\begin{document}

\maketitle

\begin{abstract}
Normalized free semi-circular random variables satisfy an upper
Khintchine inequality in $L_\infty$. We show that this implies
 the corresponding upper Khintchine inequality in any
noncommutative Banach function space. As applications, we obtain a
very simple proof of a well-known interpolation result for row and
column operator spaces and, moreover, answer an open question on noncommutative
moment inequalities concerning \cite{BeC10}.
\end{abstract}

\makeatletter
\renewcommand{\@makefntext}[1]{#1}
\makeatother

\section{Introduction}

The noncommutative analogues of Khintchine's inequalities were established by Lust-Piquard and Pisier \cite{LP86, LLP91}. Let $\M_n$ be the space of $n\times n$ matrices with complex entries and for $1\leq p\leq\infty$ let $S_p$ denote the $p$-th Schatten space. If
$\varepsilon_i$ are independent Rademacher random variables on some probability space $(\Omega,P)$, then
for any $x_i\in \M_n$ and $d\geq 1$,
$$ \Big\|\sum_{i=1}^d \varepsilon_i x_i\Big\|_{L_p(\Omega,S_p)} \approx_p
\left\{\begin{array}{ll}
\left( \Big\| \big( \sum_{i=1}^d x_ix_i^*\big)^{1/2}\Big\|_{S_p}^p+
\Big\|\big(\sum_{i=1}^d x_i^*x_i\big)^{1/2}\Big\|_{S_p}^p\right)^{\frac 1 p}&
\textrm{ if } 2\leq p<\infty \\
\displaystyle{\inf_{x_i=a_i+b_i}}\left( \Big\| \big(\sum_{i=1}^d a_ia_i^*\big)^{1/2}\Big\|_{S_p}^p+
\Big\|\big(\sum_{i=1}^d b_i^*b_i\big)^{1/2}\Big\|_{S_p}^p\right)^{\frac 1 p}&
\textrm{ if } 1\leq p<2. \end{array}\right.$$

\noindent In terms of operator spaces, these inequalities have the following interpretation. For $1\leq p \leq \infty$, let
$$Rad_p=\overline\spn \{\varepsilon_i,\, i\geq 1\}\subset L_p(\Omega)$$
and define
$$RC_p=\overline \spn\{ \delta_i=e_{i,0}+ e_{0,i};\; i\geq 1\}\subset S_p
  \ \ \mathrm{if} \ p\geq 2, \quad \ \ RC_{p} =RC_{\frac p {p-1}}^* \ \ \mathrm{if} \ 1\leq p<2,$$
where $e_{i,j}$ ($i,j\geq 0$) denote the standard matrix units. Then, for any $1\leq p<\infty$, the formal map $i_p:RC_p\to Rad_p$ defined by $i_p(\delta_i)=\varepsilon_i$ induces an isomorphism $Rad_p\approx_p RC_p$.

More generally, let $\{c_i;\,i\geq 1\}$ be a family of elements
in a finite von Neumann algebra $(\m N,\nu)$, where $\nu$ is a normal, faithful normalized trace, and set
$$G_p=\overline\spn
\{c_i;\, i\geq 1\} \subset L_p(\m N).$$

\noindent Let $1\leq p\leq \infty$. We say that the family $\{c_i;\,i\geq 1\}$
\emph{satisfies an upper Khintchine inequality in $S_p$} with constant
$C_p$ if the formal map $i_p:RC_p\to G_p$, $i_p(\delta_i)=c_i$ has
completely bounded norm smaller than $C_p$. For $2\leq p\leq\infty$, this
means that for any $x_i\in \M_n$ and $d\geq 1$
\begin{eqnarray}
\label{eqn:KIupperp}
\Big\|\sum_{i=1}^d c_i\tens x_i\Big\|_{L_p(\m N\tens \M_n)} \leq
C_p\left( \Big\| \big(\sum_{i=1}^d x_ix_i^*\big)^{1/2}\Big\|_{S_p}^p+
\Big\|\big(\sum_{i=1}^d x_i^*x_i\big)^{1/2}\Big\|_{S_p}^p\right)^{\frac 1 p},
\end{eqnarray}


\noindent with the usual modification if $p=\infty$. Typical examples of families satisfying
an upper Khintchine inequality
in $S_{\infty}$ (with constant 2) are free normalized semi-circular
random variables in the sense of Voiculescu or the family of operators
corresponding to generators of the free group $\mathbb F_\infty$ in
$L(\mathbb F_\infty)$. We refer to \cite{Pis98} for more information.

Dually, the family $\{c_i;\,i\geq 1\}$ is said to \emph{satisfy a
lower Khintchine inequality in $S_p$} with constant $C_p$ if the
formal map $j_p:G_p\to RC_p$, $j_p(c_i)=\delta_i$ has completely
bounded norm smaller than $C_p$. If the $c_i$ are orthogonal in $L_2$
and $C=\inf \{ \|c_i\|_2\}$, then $\{c_i;\,i\geq 1\}$ automatically
satisfies a lower Khintchine inequality in $S_p$ for $2\leq p \leq
\infty $:
\begin{eqnarray}
\label{eqn:KIlowerp}
C \max\left\{ \Big\|\big(\sum_{i=1}^d x_ix_i^*\big)^{1/2}\Big\|_{S_p},\,
\Big\|\big(\sum_{i=1}^d x_i^*x_i\big)^{1/2}\Big\|_{S_p} \right\}\leq
\Big\|\sum_{i=1}^d
c_i\tens x_i\Big\|_{L_p(\m N \tens\M_n)}.
\end{eqnarray}

When the $c_i$ are free semi-circular random variables, $G_p$ is
complemented in $L_p$ and the Khintchine inequalities are related to
interpolation properties of the spaces $RC_p$.

The classical approach (see \cite{Pis98}) establishes both (\ref{eqn:KIupperp}) and (\ref{eqn:KIlowerp}) at the same time
but relies on the noncommutative Khintchine inequalities for independent Rademacher variables.
In the following section, we provide an alternative very short proof. We show that
an upper Khintchine inequality in $S_{\infty}$ yields some strong estimates
in terms of distribution functions.

\section{Khintchine inequalities}

Throughout, we let $\cM$ be a von Neumann algebra equipped with a
normal, semi-finite faithful trace $\tr$. If $x$ is a closed, densely
defined operator affiliated with $\cM$, we set
$$ \lambda_t(x)=\tau (E_{(t,\infty)}(|x|)) \qquad (t\geq 0),$$ where
$E_{(t,\infty)}(|x|)$ is the spectral projection of $|x|=(x^*x)^{1/2}$
associated with $(t,\infty)$. One may think of $\lambda_t(x)$ as the
distribution function of $|x|$ with respect to $\tr$. It shares many
properties with the classical distribution functions familiar from
probability theory. Although the distribution clearly depends on the
von Neumann algebra and the trace, we do not emphasize this as the context
will always be clear. If $\lambda_t(x)<\infty$ for some $t>0$, then we
say that $x$ is $\tr$-measurable. The set $S(\tr)$ of all $\tau$-measurable
operators defines a topological $*$-algebra with respect to the
measure topology. If $x \in S(\tr)$, then the right continuous inverse
of $\lambda(x)$ given by
$$\mu_t(x) = \inf\{s>0 \ : \ \lambda_s(x) \leq t\} \qquad (t\geq 0)$$
is called the decreasing rearrangement of $x$ or the generalized
singular numbers of $x$. \par In the proof of lemma~\ref{dist} we use
the following facts. If $x\in S(\tr)$ and $t\geq 0$, then
$\lambda_t(x)=\lambda_t(x^*)$.  Moreover, for $x,\,y\in S(\tr)$, and
$s,\, t\geq 0$
$$\lambda_{t+s}(x+y)\leq \lambda_t(x)+\lambda_s(y),$$
and if $q\in \cM$ is a projection then
$$\lambda_{t}(xq)\leq \tau(q).$$
We refer to \cite{FK86} for proofs of these facts and more information.\par
Let $\cN$ denote a finite von Neumann algebra equipped
with a normal, faithful normalized trace $\nu$. We say that a sequence $\{c_i;\,i\geq 1\}\subset \cN$ 
\emph{satisfies an upper Khintchine inequality in $L_{\infty}(\cM)$} with constant $C$ if
\begin{eqnarray}
\label{eqn:KIuppervN}
\Big\|\sum_{i=1}^d c_i\tens x_i\Big\|_{\infty} \leq
C \max\left\{ \Big\| \big(\sum_{i=1}^d x_i^*x_i\big)^{1/2}\Big\|_{\infty}, 
\Big\|\big(\sum_{i=1}^d x_ix_i^*\big)^{1/2}\Big\|_{\infty}\right\},
\end{eqnarray}
for any sequence $(x_i)$ in $\cM$ and $d\geq 1$. 
\begin{lemma}\label{dist}
Suppose $\{c_i;\,i\geq 1\}\subset \cN$ satisfies an upper Khintchine inequality
in $L_{\infty}(\cM)$ with constant $C$. If $(x_i)$ is a sequence in
$S(\tr)$, then for any $t\geq 0$ and $d\geq 1$,
$$ \lambda_{Ct}\Big(\sum_{i=1}^d c_i\tens x_i \Big)\leq\lambda_t\Big(
\big(\sum_{i=1}^d x_i^*x_i\big)^{1/2}\Big)+\lambda_t\Big(
\big(\sum_{i=1}^d x_ix_i^*\big)^{1/2}\Big).$$
\end{lemma}
\begin{pf}
Fix $t\geq 0$. Let $x=\big(\sum_{i=1}^d |x_i|^2\big)^{1/2}$, $1-p=1\tens E_{(t,\infty)}(x)$, similarly $x'=\big(\sum_{i=1}^d |x_i^*|^2\big)^{1/2}$,
$1-q=1\tens E_{(t,\infty)}(x')$ and $y=\sum_{i=1}^d c_i\tens x_i$. By
the Khintchine inequality in $L_{\infty}$,
$$\| qyp\|_\infty \leq C  \max\Big\{
\| E_{[0,t]}(x) x^2 \|_\infty^{1/2},\| E_{[0,t]}(x') {x'}^{2}
\|_\infty^{1/2}\Big\}\leq  Ct.$$
We deduce that $\lambda_{Ct}(qyp)=0$. Now write
$y=qyp+ (1-q)yp +y(1-p)$ to obtain
\begin{eqnarray*}
\lambda_{Ct}(y)&\leq& \lambda_{Ct}(qyp) + \lambda_0((1-q) y p) +
\lambda_0(y(1-p))\\
& \leq & 0 + \nu\ot \tau(1-q) + \nu\ot \tau(1-p) \\
& = & \lambda_t(x)+ \lambda_t(x').
\end{eqnarray*}
\end{pf}

\noindent Let $\cM\oplus\cM$ be the direct sum von Neumann algebra
equipped with its natural trace $\tr\oplus\tr$. If $x,\,x'\in S(\tr)$,
then we can write $\lambda_t(x)+ \lambda_t(x')=\lambda_t(x\oplus x')$,
where we view $x\oplus x'$ as an element of $S(\tr\oplus\tr)$.

Let $S(0,\infty)$ be the linear space of all Lebesgue measurable, a.e.\ finite functions $f$ on $(0,\infty)$ such that $\lambda_t(f)<\infty$ for some $t>0$. Let $E$ be a symmetric space on $(0,\infty)$, i.e.\ a normed linear
subspace of $S(0,\infty)$ which is complete under its norm and
satisfies the property
$$f \in S(0,\infty), \ g \in E, \ \mu(f)\leq \mu(g) \Rightarrow f \in E \
\mathrm{and} \ \|f\|_E \leq \|g\|_E.$$
Given a symmetric space $E$ and a semi-finite von Neumann algebra, the associated noncommutative Banach
function space $E(\cM)$ is the linear subspace of all $x \in S(\tau)$
satisfying $\|x\|_{E(\cM)} := \|\mu(x)\|_E<\infty$. It is shown in
\cite{KaS08} that $E(\cM)$ is a Banach space under the norm
$\|\cdot\|_{E(\cM)}$ (see also \cite{DDdP93, Xu91} for earlier proofs
of this result under additional assumptions). In this way, we obtain a
noncommutative version
of every classical symmetric function space. In particular, if
$E=L_p$, then we recover the usual noncommutative $L_p$-space and if
moreover $\cM=B(\ell_2)$, then $E(\cM)=S_p$, the $p$-th Schatten
space.\par Taking the right-continuous inverse in lemma~\ref{dist}
yields the following result.
\begin{cor}
\label{cor:inftyKI}
Let $E$ be a symmetric space. If $\{c_i;\,i\geq 1\}$ satisfies an upper
Khintchine inequality in $L_{\infty}(\cM)$ with constant $C$, then it also
satisfies it in $E(\cM)$ with constant $C$, i.e.\ for any sequence $(x_i)$
in $E(\cM)$ we have
\begin{equation*}
\Big\|\sum_{i=1}^d c_i\ot x_i\Big\|_{E(\cN\ot\cM)} \leq C \Big\|
\big(\sum_{i=1}^d |x_i|^2\big)^{\frac{1}{2}}\oplus \big(\sum_{i=1}^d
|x_i^*|^2\big)^{\frac{1}{2}}\Big\|_{E(\cM\oplus \cM)}.
\end{equation*}
\end{cor}
Note also that, since $\lambda_t(x\oplus x')\leq \lambda_t(x + x')$
for all $x,\,x' \in S(\tr)$ and $t\geq 0$, we have
$$\Big\|
\big(\sum_{i=1}^d |x_i|^2\big)^{\frac{1}{2}}\oplus \big(\sum_{i=1}^d |x_i^*|^2\big)^{\frac{1}{2}}\Big\|_{E(\cM\oplus \cM)}\leq
\Big\|
\big(\sum_{i=1}^d |x_i|^2\big)^{\frac{1}{2}}\Big\|_{E(\cM)}+\Big\| \big(\sum_{i=1}^d |x_i^*|^2\big)^{\frac{1}{2}}\Big\|_{E(\cM)}.$$
We interpret the result in corollary~\ref{cor:inftyKI} for $E=L_p$ in terms of operator spaces.
\begin{cor}
\label{cor:interpolate}
Assume $\{c_i;\,i\geq 1\}$ is orthonormal in $L_2(\m N)$ and satisfies
an upper Khintchine inequality in $S_\infty$ with constant $C$.
For any $1\leq p\leq \infty$, $G_p\approx RC_p$ with cb constant at
most $2^{\frac 1{\max\{p,p'\}}} C$. \break \indent The family of
spaces $(RC_p)_{1\leq p\leq \infty}$ is a complex interpolation scale
(with constant at most $4\sqrt 2$).
\end{cor}
\begin{pf}
This is standard, so we only give a sketch.  For $p\geq 2$, $i_p$ has
 completely bounded norm smaller than $C$ by
 corollary~\ref{cor:inftyKI} (applied for $E=L_p$ and $\cM = \M_n$). On the other hand, the map $P_p:L_p(\m
 N)\to RC_p$ given by $P_p(x)=\sum_i \tau(c_i^* x)\delta_i$ is
 completely bounded by $2^{1/p}$. Since $i_pP_p=\Id$ on $G_p$ and
 $P_pi_p=\Id$ on $RC_p$, the result follows. This argument also shows
 that $G_p$ is completely complemented in $L_p(\m N)$ with constant
 $2^{1/p}C$. Taking adjoints gives the result for $1\leq
 p<2$. Choosing the generators of the free group as $c_i$, the
 statement about interpolation follows from the fact that $(L_p(\m
 N))_{1\leq p\leq \infty}$ is an exact complex interpolation scale (see e.g.\ \cite{PiX03}, section 2).
\end{pf}

\smallskip

\begin{rk} The interpolation result in corollary~\ref{cor:interpolate} can be extended to more general interpolation couples of noncommutative Banach function spaces, see \cite{Dir11}.\end{rk}

We obtain a very simple proof of a special case of \cite{DDPS}, theorem 4.1, with an improved constant.
\begin{cor}
Let $1\leq q<\infty$ and let $E$ be a symmetric space. If $E$ is an exact interpolation space for the couple $(L_1,L_q)$, then
\begin{equation*}
\Big\|\sum_{i=1}^d \varepsilon_i\ot x_i\Big\|_{E(L_{\infty}\ot\cM)}
\leq 4\sqrt{q} \,\Big\|\big(\sum_{i=1}^d |x_i|^2\big)^{\frac{1}{2}}\oplus \big(\sum_{i=1}^d |x_i^*|^2\big)^{\frac{1}{2}}\Big\|_{E(\cM\oplus \cM)},
\end{equation*}
for any sequence $(x_i)$ in $E(\cM)$.
\end{cor}
\begin{pf}
Let $\lambda(g_i)$ be the operators in the free group von Neumann algebra $L(\mathbb{F}_{\infty})$ corresponding to the generators of $\mathbb{F}_{\infty}$. In the proofs of \cite{LeS08}, lemma 5.4 and 5.5, it is shown using interpolation that
$$\Big\|\sum_{i=1}^d \varepsilon_i\ot x_i\Big\|_{E(L_{\infty}\ot\cM)} \leq 2\sqrt{q} \Big\|\sum_{i=1}^d \lambda(g_i)\ot x_i\Big\|_{E(L(\mathbb{F}_{\infty})\ot\cM)},$$
so the result follows immediately from corollary~\ref{cor:inftyKI}.
\end{pf}

\section{Application to moment inequalities}

We now focus on noncommutative moment inequalities associated with Orlicz functions, which were considered earlier in \cite{BeC10,BCL11}. We refer to \cite{BeC10} for the terminology used here. For an Orlicz function $\Phi$, we will need the indices
$$M(t,\Phi)=\sup_{s>0}\frac {\Phi(ts)}{\Phi(s)},\qquad
p_\Phi=\lim_{t\downarrow 0} \frac {\log(M(t,\Phi))}{\log t},\qquad
q_\Phi=\lim_{t\to \infty} \frac {\log(M(t,\Phi))}{\log t}.$$
We will assume throughout that $\Phi$ satisfies the global $\Delta_2$-condition, i.e., for some constant $C>0$,
\begin{equation}
\label{eqn:D2global}
\Phi(2t) \leq C \Phi(t) \qquad (t\geq 0).
\end{equation}

Under this condition, we automatically have $q_{\Phi}<\infty$ (\cite{KrR61}, theorem 4.1). For a semi-finite von Neumann algebra $\cM$ we let $L_{\Phi}(\cM)$ denote the
associated noncommutative Orlicz space.

From now on, we let $\{s_i;\,i\geq 1\}$ denote a family of free normalized
semi-circular variables in the finite von Neumann algebra $(\Gamma_0,\tau)$
they generate. This family satisfies (\ref{eqn:KIuppervN}) for any semi-finite 
von Neumann algebra $\cM$ \cite{Pis98}.
\begin{lemma}
\label{mom}

Let $\Phi$ be an Orlicz function satisfying (\ref{eqn:D2global}). Then, 
there is a constant $C_\Phi>0$ depending
only on $\Phi$, such that for any sequence $(x_i)$ in $L_{\Phi}(\cM)$,
$$\tau \Big(\Phi\Big(\Big|\sum_{i=1}^d s_i\tens x_i\Big|
\Big)\Big)\leq C_\Phi \max\left\{ \tau \Big(\Phi\Big(\big(\sum_{i=1}^d
x_ix_i^*\big)^{1/2}\Big)\Big);\; \tau \Big(\Phi\Big(\big(\sum_{i=1}^d
x_i^*x_i\big)^{1/2}\Big)\Big)\right\}.$$
\end{lemma}
\begin{pf}
Let $y = \sum_{i=1}^d s_i\tens x_i$, $x = \big(\sum_{i=1}^d
|x_i|^2\big)^{1/2}$ and $x' = \big(\sum_{i=1}^d
|x_i^*|^2\big)^{1/2}$. By functional calculus and lemma~\ref{dist},
$$\tau\big(\Phi(|y|)\big) = \int_0^{\infty} \Phi(\mu_t(y)) \ dt \leq
\int_0^{\infty} \Phi(2\mu_t(x\oplus x')) \ dt = \tau\big(\Phi(2(x\oplus
x'))\big).$$ Since $\Phi$ satisfies the global $\Delta_2$-condition, the result
follows.
\end{pf}

We can now answer a question formulated in \cite{BCL11} (Remark 5.3)
about moment inequalities concerning the previous paper \cite{BeC10}.
\begin{cor}
Let $\Phi$ be an Orlicz function satisfying (\ref{eqn:D2global}) and $p_\Phi>1$. 
Then, there is a constant $C_\Phi>0$ depending only on
$\Phi$ such that for any sequence $(x_i)$ in $L_{\Phi}(\cM)$,
$$\mathbb{E} \tau \Big(\Phi \Big(\Big|\sum_{i=1}^d \varepsilon_i x_i
\Big|\Big)\Big)\leq C_\Phi \max\left\{ \tau \Big(\Phi\Big(\big(\sum_{i=1}^
dx_ix_i^*\big)^{1/2}\Big)\Big);\; \tau \Big(\Phi\Big(\big(\sum_{i=1}^d
x_i^*x_i\big)^{1/2}\Big)\Big)\right\}.$$
\end{cor}
\begin{pf}
For all $1\leq p<\infty$, let $P_p$ be the projection of $L_p(\Gamma_0)$ onto $G_p$.
The maps $T_p: L_p(\Gamma_0\tens \cM)\to
L_p(\Omega;L_p(\cM))$ given by the composition of the projection
$P_p\tens Id$ and the amplification of the formal identity $k_p(s_i)=\varepsilon_i$ are
well defined and bounded by the noncommutative Khintchine inequalities.
Moreover, they do not
depend on $p$ in the sense of interpolation. Since $1<p_\Phi\leq
q_\Phi<\infty$ we find by \cite{BeC10},
Theorem 2.1, that
$$\mathbb{E}\tau \Big(\Phi \Big(\Big|\sum_{i=1}^d \varepsilon_i x_i
\Big|\Big)\Big)\leq C_\Phi\,
\tau \Big(\Phi\Big(\Big|\sum_{i=1}^d s_i\tens x_i\Big|
\Big)\Big),$$
for some constant depending only on $\Phi$ (and $q_\Phi$). Lemma~\ref{mom} gives the conclusion.
\end{pf}

By combining the previous result with \cite{BeC10}, corollary 3.1, we obtain the following Burkholder-Gundy type inequalities.
\begin{cor}
If $\Phi$ is an Orlicz function satisfying $p_\Phi>1$ and (\ref{eqn:D2global}), then there is a constant $C_\Phi>0$ depending only on $\Phi$ such that for any martingale difference sequence $(x_i)$ in $L_{\Phi}(\cM)$,
$$\mathbb{E} \tau \Big(\Phi \Big(\Big|\sum_{i=1}^d x_i
\Big|\Big)\Big)\leq  C_\Phi \max\left\{ \tau \Big(\Phi\Big(\big(\sum_{i=1}^ dx_ix_i^*\big)^{1/2}\Big)\Big);\; \tau \Big(\Phi\Big(\big(\sum_{i=1}^d x_i^*x_i\big)^{1/2}\Big)\Big)\right\}.$$
\end{cor}

\begin{rk}
The results in this section all continue to hold if we replace the
noncommutative $\Phi$-moments $\tr\big(\Phi(|\cdot|)\big)$ by the `weak
$\Phi$-moments' $\sup_{t>0}\Phi(\mu_t(\cdot))$, considered in
\cite{BCL11}. The required modifications are not difficult and left to
the interested reader.
\end{rk}

\medskip

\noindent \textbf{Acknowledgements.} The authors thank Jan van Neerven, Ben de Pagter and Quanhua Xu for
some discussions about this subject.

\bibliographystyle{acm}

\bigskip

\footnotesize{\noindent Sjoerd Dirksen, Delft Institute of Applied Mathematics, Delft University of Technology,
P.O. Box 5031,
2600 GA, Delft,
The Netherlands\\S.Dirksen@tudelft.nl

\smallskip

\n \'Eric Ricard, Laboratoire de Math\'ematiques, Universit\'e de
Franche-Comt\'e, 25030 Besan\c{c}on Cedex, France\\
eric.ricard@univ-fcomte.fr}

\end{document}